\documentclass[11pt]{amsart}

\usepackage{amsmath}
\usepackage{amssymb}
\usepackage{amsthm}
\usepackage[letterpaper]{geometry}
\usepackage{graphicx}
\usepackage{color}

\title{Super-exponential distortion of subgroups of CAT($-1$) groups}

\author{Josh Barnard}
\address{Dept.\ of Mathematics\\
	University of Oklahoma\\
	Norman, OK 73019}
\email{jbarnard, nbrady, pdani @math.ou.edu}

\author{Noel Brady}
\thanks{Noel Brady acknowledges support from  NSF grant no.\ DMS-0505707.}

\author{Pallavi Dani}

\theoremstyle{plain}
\newtheorem{thm}{Theorem}[section]

\newtheorem{prop}[thm]{Proposition}
\theoremstyle{definition}
\newtheorem{defn}[thm]{Definition}

\theoremstyle{remark}

\newcommand{\comment}[1]{}
\newcommand{\ra}{\rightarrow}
\newcommand{\disto}[2]{\delta_{#1}^{#2}}
\newcommand{\G}[2]{G_{#1}}
\newcommand{\X}[2]{X_{#1}}
\newcommand{\lk}[2]{\mathrm{Lk}(#1, #2)}
\newcommand{\dist}[3]{\mathrm d_{#1}(#2, #3)}
\newcommand{\id}[1]{\mathrm{id_{#1}}}
\newcommand{\lengthH}[1]{\ell_H(#1)}
 
\begin{document}

\begin{abstract}
We construct 2-dimensional CAT(-1) groups which contain free subgroups 
with arbitrary iterated exponential distortion, and with distortion 
higher than any iterated exponential.
\end{abstract}

\maketitle

\section{Introduction}
The purpose of this note is to produce explicit examples of CAT(-1) 
groups containing free subgroups with arbitrary iterated exponential distortion, and
distortion higher than any iterated exponential. The construction parallels 
that of Mahan Mitra in \cite{mitra} but our groups are the fundamental groups of 
locally CAT($-1$) $2$-complexes.  The building blocks used in \cite{mitra} are hyperbolic $F_3 \rtimes F_3$ groups, which are not known to be CAT($0$).  Our building blocks are graphs of groups where the vertex and edge groups are all free groups of equal rank and the underlying graph is a bouquet of a finite number of circles.  We use 
the combinatorial and geometric techniques from 
Dani Wise's version of the Rips construction \cite{wise} to ensure that our 
building blocks glue together in a locally CAT($-1$) fashion. 

One of our motivations for producing these examples was that it was not immediate from the description that the examples in \cite{mitra} had the appropriate iterated exponential distortions. Mitra's examples are graphs of groups with underlying graph a segment of length $n$, where the vertex groups are hyperbolic $F_3\rtimes F_3$ groups, and each edge identifies the kernel $F_3$ in one vertex group with the second $F_3$ factor in the adjacent vertex group.

While it is easy to see that the $n$th power of a 
hyperbolic automorphism of a free group will send a generator 
of the free group into a word which grows exponentially in $n$, it appears to be 
harder to see (without using bounded cancellation properties of carefully chosen automorphisms) that a word of length $n$ in three hyperbolic 
automorphisms (and their inverses) will send a generator of the free 
group to a word which grows exponentially in $n$.
In contrast, the monomorphisms in the multiple HNN extensions in our  
construction are all defined using positive words. This makes it
easy to see that the exponential distortions compose as required.
Also, the example in \cite{mitra} with distortion higher than any iterated exponential is of the form $(F_3\rtimes F_3)\rtimes\mathbb{Z},$
with the generator of $\mathbb{Z}$ conjugating the generators of the first $F_3$ to ``sufficiently random'' words in the generators of the second $F_3$. In contrast, our group can be described explicitly, without recourse to random words, allowing for an explicit check that our group is CAT($-1$).

Recall that if $H \subset G$ is a pair of finitely generated groups with word
metrics $d_H$ and $d_G$ respectively, the \emph{distortion} of $H$ in
$G$ is given by
$$
\disto H G(n) = \mathrm{max} \{d_H(1,h)\mid h\in H \text{ with }
d_G(1,h) \leq n\}.
$$
Up to Lipschitz equivalence,
this function is independent of the choice of word metrics.   
Background on CAT($-1$) spaces and the large link condition 
may be found in \cite{BH}.

\section{The building blocks}

For each positive integer $n$ we define a building block group  
\[
\G {n}{m} = \langle a_1, \dots, a_m, t_1, \dots, t_n \mid t_i a_jt_i^{-1} =
W_{ij}\; ; \;  1 \leq i \leq n, 1\leq j\leq m  \rangle,
\]
where $m=14n$ and $\{W_{ij}\}$ is a collection of positive words of length $14$ 
in $a_1,\dots, a_n$,
such that each two-letter word $a_ka_l$ appears at most once among the $W_{ij}$'s.
One way of ensuring this is to choose these words to be consecutive subwords of the 
following word, defined by Dani Wise in \cite{wise}.

\begin{defn}[Wise's long word with no two-letter repetitions] \label{long-word}
Given the set of letters
$\{a_1, \dots, a_m\}$, define
\[
\Sigma(a_1, \dots, a_m)=(a_1a_1a_2a_1a_3\cdots a_1a_m)(a_2a_2a_3a_2a_4\cdots
a_2a_m) \cdots\cdots(a_{m-1}a_{m-1}a_m)a_m.
\]
\end{defn}

It is easy to see that $\Sigma(a_1, \dots, a_m)$ is a positive word of length 
$m^2$, such that each two-letter word $a_k a_l$ appears as a subword in at 
most one place.  Following Dani Wise, we simply chop $\Sigma(a_1, \ldots, a_m)$ 
into subwords of length 14. 
 In order to obtain all $mn$ relator words of $G_n$  from $\Sigma(a_1,\ldots,a_m)$ 
we must have $m^2 \geq 14mn$, which explains our choice of $m$ above.

\begin{prop}\label{build-block}
The presentation $2$-complex $\X n {m}$ for $\G{n}m$ can be given a locally CAT($-1$) structure.  Furthermore,
\begin{enumerate}
\item The $a_j$'s generate a free subgroup $F(a_j)$ whose distortion in $\G n m$ is exponential.
\item The $t_i$'s generate a free subgroup $F({t_i})$ of $\G n m$ that is highly convex in the following sense:

\noindent
Let $v$ be the vertex of $\X n m$.  Then 
\[
\dist{\lk{v}{\X n m}}{t_i^{\epsilon_1}}{t_j^{\epsilon_2}} \geq 2\pi, 
\text{ where } \epsilon_1, \epsilon_2 \in \{+,- \} 
\text{ and if } 
i=j 
\text{ then }  
\epsilon_1 \neq \epsilon_2.
\]
\end{enumerate}
\end{prop}

\begin{figure}[!ht]
\begin{center}
\scalebox{0.85}{\input{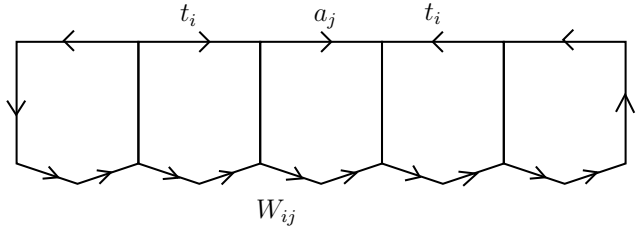}}
\caption{A 2-cell of $X_{n}$ decomposed into right-angled pentagons.
\label{2cell}}
\end{center}
\end{figure}

\begin{proof}
Each disk in $\X nm$ is given a piecewise hyperbolic structure 
by expressing it as a concatenation of right-angled hyperbolic 
pentagons, 
as shown in 
Figure \ref{2cell}.  
The fact that $\X nm$ satisfies the large link condition 
is a consequence of the condition that the $W_{ij}$'s are positive 
words with no two-letter repetitions. The details are in \cite{wise}. 
Thus $\X nm$ is locally CAT($-1$) and $G_n$ is hyperbolic. 

To prove (1), first note that the group $\G nm$ can be viewed as the fundamental group of a graph of groups. The underlying graph is a bouquet of $n$ circles, and the edge and vertex groups are all equal to $F({a_j})$, the free group on $a_1,\dots, a_m$. The two maps associated with the $i$th edge are 
$\id{F({a_j})}$ and $\phi_i: F({a_j}) \ra F({a_j})$, defined by $\phi_i(a_j)=W_{ij}$ for $1\leq j \leq m$. 
To see that $\phi_i$ is injective,  
note that $\phi_i$ induces a map on a subdivided bouquet of $m$ circles.
By Stallings' algorithm this factors through a sequence of folds followed by an immersion.  	The no two-letter repetitions 
condition restricts the amount of folding that can occur and 
ensures that no non-trivial loops are killed by the sequence of folds.  Thus the 
folding maps induce isomorphisms at the level of $\pi_1$ and
the final immersion induces an injection. 
The theory of graphs of groups now implies that
the subgroup of $\G nm$ generated by $a_1, \dots, a_m$ is 
$F({a_j})$.  
 
The distortion of $F({a_j})$ in $\G nm$ is at least exponential because, for example, the element $t_1^ka_1t_1^{-k}$ of $F({a_j})$, when expressed in terms of the $a_j$'s, is a positive (hence reduced) word of length $14^k$. 
To see that the distortion is at most exponential consider a word $w(a_j, t_i)$ which represents an element of $F({a_j})$ and has length $k$ in $\G nm$.  It can be reduced to a word in the $a_j$'s by successively cancelling at most $k/2$ innermost $t_i\cdots t_i^{-1}$ pairs. Each such cancellation multiplies the word length by at most a factor of $14$. So the length of $w(a_j, t_i)$ in $F(a_j)$ is at most $14^{k/2}k$.

To prove (2), note that any path in the link of $v$ connecting $t_i^{\pm}$ 
and $a_j^{\pm}$ within a disk as in Figure \ref{2cell}
has combinatorial length $2$. So any path in the link connecting $t_i^{\epsilon_1}$ and $t_j^{\epsilon_2}$ has combinatorial length at 
least $4$. Since we are using right-angled pentagons, combinatorial length 
4 corresponds to a spherical metric length of $4(\pi/2)$ or $2\pi$. 
Thus $\langle t_1, \dots, t_n \rangle$ is highly convex.  
As a consequence, we have that the map from the bouquet of circles 
with edges labelled $t_1, \dots, t_n$ into $\X nm$ is a 
local isometric embedding. This
implies that it is $\pi_1$-injective, i.e. $t_1, \dots, t_n$ generate the free group $F({t_i})$.  To see this algebraically, note that  
the homomorphism $\psi:\G nm \ra F({t_i})$ defined by $\psi(t_i)=t_i$ and $\psi(a_j)=1$ is a retraction of $G_n$ onto $F({t_i})$.
\comment{
This can also be seen algebraically as follows. 
There is an obvious homomorphism $\phi$ from $F({t_i})$ into $\G nm$ (namely $\phi(t_i)=t_i$ for all $i$).  Also, the map $\psi:\G nm \ra F({t_i})$ defined by $\psi(t_i)=t_i$ and $\psi(a_j)=1$ is a homomorphism such that $\psi \circ \phi =\mathrm{id_{F({t_i})}}$.  Thus $\phi$ is injective.}   
\end{proof}

\section{Iterated exponential distortion}

In this section we see how to get arbitrary iterated exponential distortions. 
The idea is to amalgamate a chain of building block groups together, 
identifying the distorted free group in one with the highly convex free 
subgroup in the next. This identification of distorted with highly convex 
can be made in a non-positively curved way, and the distortion functions 
compose as expected. Here are the details.

\begin{thm}
For any integer $l>0$, there exists a $2$-dimensional CAT($-1$) group $H_l$
with a free subgroup $F$ such that $\disto{F}{H_l}(x) \simeq \exp^l
(x)$.
\end{thm}
We will actually show that $\disto{F}{H_l}(x) \simeq f^l(x)$, where
$f(x)=14^x$.

\begin{proof}
The group $H_l$ is defined using the building blocks from Proposition \ref{build-block} as

\[
H_l = \G{1}{14} \ast_{F_{1}} \G{14}{14^2} \ast_{F_{2}} \G{14^2}{14^3} 
\ast  \cdots \cdots \ast_{F_{{l-1}}} \G{14^{l-1}}{14^l},
\]
where $F_k$, for $1 \leq k\leq  l-1$, is a free group of rank $14^k$ which is 
identified with the exponentially distorted free subgroup of 
$\G{14^{k-1}}{14^{k}}$  and the highly convex free subgroup of 
$\G{14^{k}}{14^{k+1}}$. 
Let $a_j^{(k)}$, with $1 \leq j \leq 14^k$, denote the generators of $F_k$, and 
let $t$ denote the stable letter of $G_1$, 
so that $G_1 = \langle a_j^{(1)}, t \rangle$ and
$G_{14^k}= \langle a_j^{(k+1)}, a_j^{(k)} \rangle$. We shall use this notation 
in the upper and lower bound arguments below. 

Let $Y_l$ denote the presentation complex of $H_l$.  Then $Y_1=X_1$, which is locally CAT($-1$) by Proposition \ref{build-block}.  Further, there are inclusions $Y_1 \subset Y_2 \subset Y_3 \cdots$, so that 
the large link condition can be checked inductively.  The space $Y_{k+1}$ 
is obtained by gluing $X_{14^{k}}$ to $Y_{k}$ along a rose $R_{k}$ 
with $14^{k}$ petals, and $R_k$ is highly convex in $X_{14^k}$.  It follows that
the link of the base vertex of $Y_{k+1}$ is obtained by gluing together the links of 
base vertices  in $Y_{k}$ and $X_{14^{k}}$, along a set of $2(14^k)$ points which is 
$2\pi$-separated in the latter link. By induction, the link of the base vertex in 
$Y_k$ is large, and hence the union is large. 

The group with the desired distortion in $H_l$ is $F_l$, the free group generated by $a_j^{(l)}$, with $1\leq j \leq 14^l$, which is exponentially distorted in $G_{14^{l-1}}$.  
We prove the lower bound as follows.  Given a positive integer  $n$
consider the sequence given by
$w_1 = t^na_1^{(1)}t^{-n} \in F_1$ and $w_{k}=w_{(k-1)}a_1^{(k)}w_{(k-1)}^{-1}\in F_k$, and set $w=w_l$.
Given $g \in H_l$ let $\ell_{H_l}(g)$ denote the distance from $1$ to $g$ in $H_l$.  Observe that $\ell_{H_l}(w_1)\leq 2n+1$, and inductively that $\ell_{H_l}(w)\leq 2^ln+2^l-1$ (a linear function of $n$). 
On the other hand, each $w_k$ can be expressed as a positive word in the generators of $F_k$ by using the $W_{ij}$'s from the definition of the building blocks.  Since there is no cancellation among positive words, we obtain that $|w_1|_{F_1}=14^n$ and inductively that  $|w|_{F_l}=f^l(n)$. 
This gives the lower bound $f^l(x)\preceq \disto{F}{H_l}(x) $.

To prove the upper bound it is more convenient to do the induction in the opposite direction.  Proposition \ref{build-block} provides the base case. 
Let $w$ be an element whose length in $H_l$ is at most $n$.  Then by successively cancelling at most $n/2$ innermost $t\cdots t^{-1}$ pairs, $w$ can be 
represented by a word in $\G{14}{14^2}\ast_{F_{2}}G_{14^2} \ast \cdots \G{14^{l-1}}{14^l}$ of length at most $14^{n/2}n$.  
At each stage of the previous cancellation procedure, we may assume (by 
replacement if necessary) that 
the subword enclosed by an innermost $t\ldots t^{-1}$ involves only the $a^{(1)}_i$. 
This is because the free group on the $a_i^{(1)}$ is a retract of 
  $\G{14}{14^2}\ast_{F_{2}}G_{14^2} \ast \cdots \G{14^{l-1}}{14^l}$. 
  The upper bound follows by induction.  
\end{proof}

\section{Distortion higher than any iterated exponential}

In this section we produce a 2-dimensional CAT(-1) group containing a 
free subgroup with distortion more than any iterated exponential. The idea is to 
take a suitable building block group from section 2 with base group free on 
$a_i$ and stable letters $t_j$, and to form a new HNN extension with stable letter 
$s$ which sends the free group on the $a_i$ into  a subgroup of the free 
group on the $t_j$. 
 
\begin{thm}
There exists a $2$-dimensional CAT($-1$) group $G$ with a subgroup $H$, such that $\disto H G (x)$ is a function that is bigger than any iterated
exponential.
\end{thm}

\begin{figure}[!ht]
\begin{center}
\scalebox{0.85}{\input{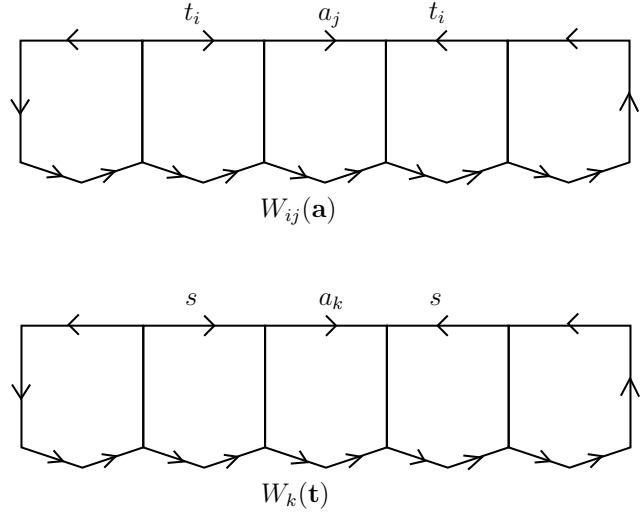}}
\caption{Relator 2-cells of the group $G$ decomposed into right-angled
pentagons. \label{second}}
\end{center}
\end{figure}

\begin{proof}
Define
\[
G=\left \langle a_1, \dots, a_m, t_1,\dots, t_n, s \mid
t_ia_jt_i^{-1}=W_{ij}\; ;\; sa_ks^{-1}=W_k \right\rangle,
\]
where $\{W_{ij}\}$ (resp. $\{W_k\}$) consists of $mn$ (resp. $m$) positive
words of length $14$ in the letters $a_i$ (resp. $t_j$), with no
$2$-letter repetitions.  Thus we may
choose the $W_{ij}$'s and $W_k$'s to be disjoint subwords of
$\Sigma(a_1, \dots, a_m)$ and $\Sigma(t_1, \dots, t_n)$ (see Definition \ref{long-word}) respectively.
This gives the following two conditions:
\[
14mn \leq m^2 \quad \text{ and } \quad 14m \leq n^2.
\]
So, for example, $n=14^2$ and $m=14^3$ is a possible choice.

\begin{figure}[!ht]
\begin{center}
\scalebox{0.6}{\input{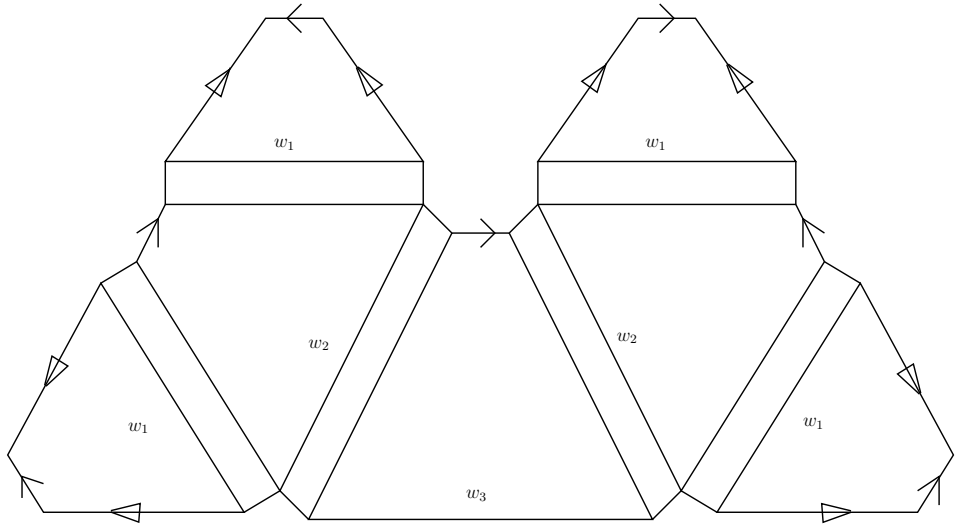}}
\caption{ The word $w_k$ (in this case $k=3$). \label{snow}}
\end{center}
\end{figure}

Each disk in the presentation complex is given a piecewise hyperbolic structure  by 
concatenating a sequence of right-angled pentagons as  shown  
in Figure \ref{second}. Observe that loops in the link of the unique vertex
which involve $s^{\pm}$ have combinatorial length at least $4$.  Similarly, any loop  involving an $a_j^{\pm}$
and a $t_i^{\pm}$, for some $i$ and $j$, has length at least 4.  Thus to show that the
complex is locally CAT($-1$), it is enough to consider loops involving only $a_i^{\pm}$'s or only $t_j^{\pm}$'s. The fact that such loops are large is a consequence of the condition that the $W_{ij}$'s and $W_k$'s are positive words with no two-letter repetitions. The details are exactly as in \cite{wise}.

Let $H$ be the subgroup of $G$ generated by $a_1, \dots, a_m$.  Then
$\disto H G (x)$ is higher than any iterated exponential.  To see
this, consider the sequence of words $w_k \in H$ given by
$w_1 =t_1a_1t_1^{-1}$ and $w_{k}=(sw_{k-1}s^{-1}) a_1 (s
w_{k-1}^{-1}s^{-1})$ for $k >1$. 
The word $w_3$ is shown in Figure~\ref{snow}. The label of each single arrow edge 
is $a_1$, the label of each solid arrow edge is $t_1$, and the 
edges along the strips are all labeled by $s$ and are oriented from 
$w_j$ toward $w_{j-1}$.  
Let $\ell_G$ and $\ell_H$ denote geodesic lengths in $G$ and $H$ respectively. 
Note that $\ell_G(w_1)=3$ and $\ell_G(w_k) \leq 2 \ell_G(w_{k-1})+5$ for $k>1$.
So, for example, $\ell_G(w_k) \leq 4^k$ is true. 
On the other hand, $\lengthH{w_1} =14$ and $\lengthH{w_k}= 14^{14\lengthH{w_{k-1}}} >14^{\lengthH{w_{k-1}}}$. So by induction $\lengthH{w_k} \geq f^k(1)$.  (Recall that $f(x)=14^x$.)  This shows that 
$\disto H G (x)\geq f^{\lfloor\mathrm{log}_4 x\rfloor}(1)$, which is a function that grows faster than 
any iterated exponential.  
\end{proof}


\begin{thebibliography}{WW}

\bibitem[BH]{BH}
M. R. Bridson and A. Haefliger, {\it Metric spaces of non-positive curvature}. Grundlehren der Mathematischen Wissenschaften, 319. Springer-Verlag, Berlin, 1999.



\bibitem[Mi]{mitra}
M. Mitra, {\it Cannon-Thurston maps for trees of hyperbolic metric spaces}, J. Diff. Geom. 48 (1998),  no. 1, 135--164.

\bibitem[Wi]{wise}
D. Wise, {\it Incoherent negatively curved groups}, Proc. Amer. Math. Soc.  126  (1998),  no. 4, 957--964.

\end{thebibliography}
\end{document}